\documentclass[12pt]{amsart}
\usepackage[cmtip,arrow]{xy}
\usepackage{
amsxtra,amsmath,amsthm,amssymb,amscd,ascmac,amsfonts,
multicol,tabularx,latexsym,mathrsfs}
\usepackage{amsfonts}
\theoremstyle{definition}

\newcommand{\Q}{\mathbb Q}
\newcommand{\Z}{\mathbb Z}

\newcommand{\Gal}{\mathrm{Gal}}
\newtheorem{thm}{Theorem}

\newtheorem{lem}{Lemma}
\newtheorem{prop}{Proposition}

\newtheorem{cor}{Corollary}

\newcommand{\Qb}{\overline{\mathbb Q}}

\newcommand{\Tor}{\mathrm{Tor}}

\newcommand{\F}{\mathbb F}

\newdimen\minCDarrowwidth
\date{}
\begin{document}
\title[]
{A number field analogue of Weil's theorem on congruent zeta functions}
\author{Manabu\ OZAKI}
\maketitle
\section{Introduction}
Let $K$ be a function field of one variable over a finite field
of characteristic $l>0$ with the constant field $\F:=K\cap\overline{\F_l}$.
We denote by $\zeta_{K}(s)$ the congruent (Dedekind) zeta function of $K$.
Furthermore let $X_{\overline{\F}K}(p)$ be the Galois group
of the maximal unramified abelian $p$-extension over $\overline{\F}K$, $p\ne l$ being any prime number.
Then Weil's celebrated theorem
states that $\zeta_{K}(s)$ is completely described by
$\Gal(\overline{\F}/\F)(\simeq\Gal(\overline{\F}K/K))$-module structure of $X_{\overline{\F}K}(p)\otimes_{\Z_p}\Q_p$.
Specifically, we have
\[
\zeta_{K}(s)=\frac{P(q^{-s})}{(1-q^{-s})(1-q^{1-s})}\ \ (s\in\mathbb{C}),
\]
where $q=\#\F$ and 
\[
P(T)=\det(1-T\,\mathrm{Frob}_{\F}|X_{\overline{\F}K}(p)\otimes_{\Z_p}\Q_p)\in\Z[T],
\]
$\mathrm{Frob}_{\F}\in\Gal(\overline{\F}/\F)$ being the Frobenius automorphism over $\F$.

\

The above theorem means that the
congruent zeta function of $K$ is determined by the $\Q_p[\Gal(\overline{\F}/\F)]$-module structure of 
$X_{\overline{\F}K}(p) \otimes_{\Z_p}\Q_p$ and
vise versa since $X_{\overline{\F}K}(p)\otimes_{\Z_p}\Q_p$ is a semisimple $\Q_p[\Gal(\overline{\F}/\F)]$-module.

In the present paper, we will give a number field analogue of Weil's theorem
considering the {\it total cyclotomic extension}, which is a number field analogue of the constant field extension $\overline{\F}K$ of a function field $K$ over $\F$.
 
Let $k$ be a number field of finite degree, and
we put $\tilde{k}:=k(\mu_\infty)$, where $\mu_\infty$ stands for
all the roots of unity in a fixed algebraic closure $\Qb$ of $\Q$.
We denote by $X_F(p)$ the Galois group of the maximal unramified
abelian $p$-extension over $F$ for any number field $F$ and prime number $p$.
Our main theorem is the following:
\begin{thm}
Let $p$ be a prime number, and let 
$k_1$ and $k_2$ be number fields of finite degree.
We assume that $F:=k_1\cap\tilde{\Q}=k_2\cap\tilde{\Q}$ and that
$N\cap\Q(\mu_p)=\Q$, where $N$ is the Galois closure of $k_1k_2$ over $\Q$.
If $X_{\tilde{k_1}}(p)\simeq X_{\tilde{k_2}}(p)$ as $\Z_p[[\Gal(\tilde{\Q}/F)]]$-modules,
then $\zeta_{k_1}(s)=\zeta_{k_2}(s)$ holds.
Here we regard $X_{\tilde{k_i}}(p)$ as $\Z_p[[\Gal(\tilde{\Q}/F)]]$-module via the isomorphism
$\Gal(\tilde{k_i}/k_i)\simeq\Gal(\tilde{\Q}/F)$ given by the restriction
($i=1,2$).
\end{thm}
We will also give a partial result for the converse of Theorem 1:
\begin{thm}
Let $k_1$ and $k_2$ are number fields of finite degree, and
let $N$ be the Galois closure of $k_1k_2$ over $\Q$. Let $p$ be a prime number.
Assume that $N\cap\tilde{\Q}=\Q$ and $p\nmid [N:\Q]$.
Then $\zeta_{k_1}(s)=\zeta_{k_2}(s)$ implies
$X_{\tilde{k_1}}(p)\simeq X_{\tilde{k_2}}(p)$ as 
$\Z_p[[\Gal(\tilde{\Q}/\Q)]]$-modules.
\end{thm}
Though the above theorems provide no explicit relationship
between the $\Z_p[[\Gal(\tilde{k}/k)]]$-module structure of $X_{\tilde{k}}(p)$ and the Dedekind zeta function $\zeta_k(s)$, we may regard it as a considerable step toward number field analogue of Weil's theorem.
\section{Proof of Theorem 1}
In this section, we will give a proof of Theorem 1.
\begin{prop}
Let $k$ be a number field of finite degree and 
$L$ a intermediate field of $\tilde{k}/k$. Assume that
for each intermediate field $M$ of $L/k$ with $[M:k]<\infty$ there exist
infinitely many degree one primes $\mathfrak{l}$ of $M$ such that
$\mu_l\subseteq L$, $l$ being the rational prime below $\mathfrak{l}$.
Then we have 
\[
\varprojlim_{k\subseteq M\subseteq L,\,[M:k]<\infty} E_M=0,
\]
where $E_M$ denotes the global unit group of $M$ and
the projective limit is taken with respect to the norm maps.
\end{prop}
{\bf Proof.}\ \ \ 
It is enough to show that 
\[
\bigcap_{M\subseteq N\subseteq L,\ [N:M]<\infty} N_{N/M}(E_N)=1
\]
for each intermediate field $M$ of $L/k$
with $[M:k]<\infty$.
Let $\frak{l}$ be a degree one prime of $M$ such that
$\mu_l\subseteq L$, $l$ being the rational prime below $\mathfrak{l}$,
and $\frak{l}$ is unramified in $M/\Q$.
Then $M(\mu_l)\subseteq L$ and 
$N_{M(\mu_l)/M}(\varepsilon)\equiv 1\pmod{\frak{l}}$ holds
for any $\varepsilon\in E_{M(\mu_l)}$.
Since $M(\mu_l)\subseteq L$, we find that
$\eta\equiv 1\pmod{\frak{l}}$ holds for any
$\eta\in\bigcap_{M\subseteq N\subseteq L,\ [N:M]<\infty}
N_{N/M}(E_N)$.
Because there exists infinitely many such primes $\frak{l}$,
we conclude that $\eta=1$ must hold.
\hfill$\Box$

\

Let $U_{M,\frak{l}}(p)$ be the pro-$p$-part of the local unit group at the prime $\frak{l}\nmid p$
of a number field $M$ of finite degree. Here we note that
$U_{M,\frak{l}}(p)$ is a finite cyclic $p$-group.
Define
\[
\mathcal{U}_{L,l}(p):=\varprojlim_{M\subseteq L, [M:\Q]<\infty}
\prod_{\frak{l}\in S_l(M)}U_{M,\frak{l}}(p)
\]
to be the projective limit of the pro-$p$-part of the semi-local unit groups of $M$ at $l$, $l$ being a prime number with $l\ne p$,
with respect to the norm maps, 
where $S_l(M)$ stands for the set of all the primes of $M$ lying over $l$
and $M$ runs over
all the subfields of $L$ with $[M:\Q]<\infty$.
\begin{prop}
Let $L\subseteq\tilde{k}$ be as in Proposition 1.
Let $l$ and $p$ be distinct prime numbers, and 
denote by $L_{\{l\}}(p)/L$ and $L_\emptyset(p)/L$ the maximal $l$-ramified and unramified abelian $p$-extension, respectively.
Then the reciprocity map induces the isomorphism
\[
\mathcal{U}_{L,l}(p)\simeq
\Gal(L_{\{l\}}(p)/L_\emptyset(p)).
\]
\end{prop}
{\bf Proof.}\ \ \ We get the following exact sequence by using class field theory:
\[
\varprojlim_{M\subseteq L,\,[M:\Q]<\infty}E_M\longrightarrow
\mathcal{U}_{L,l}(p)\longrightarrow
\Gal(L_{\{l\}}(p)/L_\emptyset(p))\longrightarrow 0.
\]
It follows from Proposition 1 that the first term of the above
exact sequence vanishes, which implies the Proposition.
\hfill$\Box$

\

Let $p$ and $l$ be prime numbers with $l\equiv 1\pmod{p}$, and
define $L_0$ to be the subfield of $\Q(\mu_l)$ with $[\Q(\mu_l):L_0]=p$.
Let $\tilde{k}^{(l')}$ be the extension of $k$ adjoining all the roots of unity whose order is prime to $l$.
We denote by $L$ the composite of $\tilde{k}^{(l')}$, the cyclotomic $\Z_l$-extension of $\Q$, and $L_0$.
We find that if $l$ is unramified in $k/\Q$ then
$\tilde{k}/L$ is a cyclic extension of degree $p$ and the 
ramified primes in this extension are exactly all the primes lying over $l$. Furthermore, $L$ satisfies the assumption of Propositions 1 and 2.
Let $\Delta:=\Gal(\tilde{k}/L)$
and put $X:=X_{\tilde{k}}(p)_\Delta=\Gal(Z/\tilde{k})$,
where $Z$ is the central $p$-class field of $\tilde{k}/L$.
Put $Y:=\Gal(Z/L)$. Then $Y$ is an abelian pro-$p$-group, on which 
$G:=\Gal(L/k)$ acts.
\begin{prop}
Assume that a prime $l$ is unramified in $k/\Q$ and
\[
l\mathcal{O}_k=\frak{l_1}\frak{l}_2\cdots\frak{l}_g
\]
is the
prime decomposition of the principal ideal $l\mathcal{O}_k$
in $k$.
Denote by $D_i$ the decomposition subgroup of $G$ for the prime $\frak{l}_i\ 
(1\le i\le g)$.

Let $\mathcal{I}\subseteq Y$ be the subgroup generated by the inertia subgroups of $Y$ for the primes of $L$ lying over $l$.
Then
\[
\mathcal{I}\simeq\bigoplus_{i=1}^g\F_p[[G/D_i]]
\]
as $\F_p[[G]]$-modules and
$\mathcal{I}=\mathrm{Tor}_{\Z_p}(Y)$.
\end{prop}
{\bf Proof.}\ \ \ 
We first note that the ramified primes of $L$ in the abelian pro-$p$ extension $Z/L$ are exactly
all the primes lying over $l$, and whose ramification indexes are $p$.
Then it follows from Proposition 2 that
\begin{equation}\label{I}
\mathcal{I}\simeq\mathcal{U}_{L,l}(p)/p,
\end{equation}
since $L_\emptyset(p)\subseteq Z$.
If we put 
\[
\mathcal{U}_{L,\frak{l}_i}(p):=\varprojlim_{k\subseteq M\subseteq L, [M:k]<\infty}
\prod_{\frak{l}\in S_{\frak{l}_i}(M)}U_{M,\frak{l}}(p),
\]
where $S_{\frak{l}_i}(M)$ is the set of all the primes of $M$
lying over $\frak{l}_i$ and the projective limit is taken with respect to the norm maps,
then we have
\begin{equation}\label{U_l}
\mathcal{U}_{L,l}(p)=\bigoplus_{i=1}^g\mathcal{U}_{L,\frak{l}_i}(p).
\end{equation}

Since $\mu_p\subseteq\Z_l^\times$, 
$U_{M,\frak{l}}(p)/p\simeq\F_p$ as $D_i$-modules
for any finite extension $M/k$ with $M\subseteq L$
and prime $\frak{l}|\frak{l}_i$ of $M$.
Hence we have
\begin{equation}\label{U}
\bigoplus_{\frak{l}|\frak{l}_i}U_{M,\frak{l}}(p)/p
\simeq
\F_p[\Gal(M/k)/D_i(M/k)]
\end{equation}
as $\F_p[\Gal(M/k)]$-modules, where $D_i(M/k)$ stands for the decomposition subgroup of $\Gal(M/k)$ for the prime $\frak{l}_i$.
Furthermore, if a finite subextension $M/k$ of $L/k$ is large enough so that the primes of $M$ over $l$ are not ramified in any finite subextension of $L/M$ with the ramification index divisible by $p$,
the local norm map $U_{M'',\frak{l}''}(p)/p\longrightarrow
U_{M',\frak{l}'}(p)/p$ is an isomorphism for any finite extensions $M'$ and $M''$ of $M$ with $M'\subseteq M''$ and primes $\frak{l}''|\frak{l}'|l$
of $M''$ and $M'$, respectively.
Hence we find that
\[
\mathcal{U}_{L,\frak{l}_i}(p)/p\simeq\F_p[[G/D_i]].
\]
Thus we conclude by using \eqref{I}, \eqref{U_l} and \eqref{U} that
\[
\mathcal{I}\simeq\bigoplus_{i=1}^g\F_p[[G/D_i]]
\]
as $\F_p[[G]]$-modules.

We have the exact sequence
\[
0\longrightarrow \mathcal{I}\longrightarrow Y\longrightarrow
\Gal(L_{\emptyset}(p)/L)\longrightarrow 0.
\]
Since $\Gal(L_{\emptyset}(p)/L)$ is a free pro-$p$ abelian group
by \cite[Theorem 1]{U} 
(because $L$ contains $\Q^{(p)}$ defined in \cite{U}),
the above exact sequence of pro-$p$ abelian groups
splits and we conclude that $\mathcal{I}=\mathrm{Tor}_{\Z_p}(Y)$.
\hfill$\Box$

\

In what follows we will derive the number of the primes of $k$
lying over $l$ from the $\Z_p[[\Gal(\tilde{k}/k)]]$-module structure of $X_{\tilde{k}}(p)$.

Put $\mathcal{J}:=X\cap\mathcal{I}$.
Then $\mathcal{J}=\mathrm{Tor}_{\Z_p}X$ by Proposition 3 and we get the exact sequence
\begin{equation}\label{ex}
0\longrightarrow \mathcal{J}\longrightarrow\mathcal{I}\longrightarrow\F_p
\longrightarrow 0
\end{equation}
of $\F_p[[G]]$-modules.

Let $I_G$ be the augmentation ideal of $\F_p[[G]]$.
We put $N_{\sigma D_i}:=\sum_{t=0}^{f-1} \sigma^tD_i\in\F_p[[G/D_i]]$
when $\sigma D_i\in G/D_i$ is of finite order $f$.

\begin{lem}
(1)\ \ \ Let $\sigma\in G$ so that $\sigma D_i\in G/D_i$ is of finite order.
Then 
\[
\F_p[[G/D_i]]^{\langle\sigma\rangle}=N_{\sigma D_i}\F_p[[G/D_i]].
\]
Furthermore, we have
\[
\ker(\F_p[[G/D_i]]\overset{N_{\sigma D_i}}{\longrightarrow}\F_p[[G/D_i]])
=(\sigma-1)\F_p[[G/D_i]].
\]
(2)\ \ \ If $\sigma\in G$ satisfies that the order of $\sigma D_i\in G/D_i$
is finite and divisible by $p$, then
\[
\F_p[[G/D_i]]^{\langle\sigma\rangle}\subseteq I_G\F_p[[G/D_i]].
\]
(3)\ \ \ If the order of $\sigma D_i$ is finite, then
\[
\F_p[[G/D_i]]^{\langle\sigma\rangle}/I_G(\F_p[[G/D_i]]^{\langle\sigma\rangle})
\simeq\F_p.
\]
\end{lem}
{\bf Proof.}\ \ \
(1)\ \ \ Since $G/D_i$ is a pro-finite group, we assume that $G/D_i=\varprojlim_{s\in S}G_s$ for a certain projective system
$\{G_s\}_{s\in S}$ of finite groups.
We note that $\langle\sigma D_j\rangle\hookrightarrow G_s$
if $s$ is large enough because $\langle\sigma D_j\rangle$ is finite.
Hence $\F_p[G_s]^{\langle \sigma\rangle}=N_{\sigma D_j}\F_p[G_s]$ for all large enough $s\in S$, which comes from $\hat{H}^0(\langle\sigma D_i\rangle,\,\F_p[G_s])=0$, implies the first assertion.

Similarly the second assertion follows from $\hat{H}^{-1}(\langle\sigma D_i\rangle,\,\F_p[G_j])=0$.

(2)\ \ \ The assertion follows from (1) and the fact that $N_{\sigma D_i}\in
I_{G/D_j}$, where $I_{G/D_j}=I_G\F_p[[G/D_i]]$
is the augmentation ideal of $\F_p[[G/D_i]]$,
if the order of $\sigma D_i\in G/D_i$ is finite and divisible by $p$.

(3)\ \ \ Since there exists a surjection 
\begin{align*}
\F_p\simeq
\F_p[[G/D_i]]/I_G\F_p[[G/D_j]]&\longrightarrow
N_{\sigma D_j}\F_p[[G/D_i]]/I_GN_{\sigma D_i}\F_p[[G/D_i]]\\
&=\F_p[[G/D_i]]^{\langle\sigma\rangle}/I_G(\F_p[[G/D_i]]^{\langle\sigma\rangle}),
\end{align*}
it is enough to show that $\F_p[[G/D_i]]^{\langle\sigma\rangle}/I_G(\F_p[[G/D_i]]^{\langle\sigma\rangle})\ne 0$.

Suppose that 
$N_{\sigma D_i}\F_p[[G/D_i]]/I_GN_{\sigma D_i}\F_p[[G/D_i]]=0$.
Then $N_{\sigma D_i}\in I_GN_{\sigma D_i}\F_p[[G/D_i]]=I_{G/D_i}N_{\sigma D_i}$.
Hence $N_{\sigma D_i}=\alpha N_{\sigma D_i}$ holds for a certain 
$\alpha\in I_{G/D_i}$, which in turn implies $(\alpha-1)N_{\sigma D_i}=0$.
However, since $\alpha-1\not\in I_{G/D_i}$, this contradicts to the second assertion of (1). Therefore  
\begin{align*}
\F_p[[G/D_i]]^{\langle\sigma\rangle}/I_G(\F_p[[G/D_i]]^{\langle\sigma\rangle})&=N_{\sigma D_j}\F_p[[G/D_i]]/I_GN_{\sigma D_i}\F_p[[G/D_i]]\\
&\ne 0.
\end{align*}
\hfill$\Box$

\begin{prop}
Let $k_1$ and $k_2$ be number fields of finite degree which satisfy
the assumption of Theorem 1.
Let $l$ be a prime number such that $l\equiv 1\pmod {p}$ and
$l$ is unramified in $k_1k_2/\Q$.
Then the numbers $g_{k_1}(l)$ and $g_{k_2}(l)$ of the primes of $k_1$ and $k_2$ lying over $l$, respectively,
coincide.
\end{prop}
{\bf Proof.}\ \ \ 
Recall that $F:=k_1\cap\tilde{\Q}=k_2\cap\tilde{\Q}$, and
that we regards $X_{\tilde{k_i}}(p)$ as 
$\mathcal{G}:=\Gal(\tilde{\Q}/F)$-modules via the isomorphism 
$\pi_i:\Gal(\tilde{k_i}/k_i)\simeq\mathcal{G}$ induced by the restriction
($i=1,2$).

Put $M:=\Q(\bigcup_{l\nmid m}\mu_m)\Q_\infty^{(l)}$, 
$\Q_\infty^{(l)}/\Q$ being the cyclotomic $\Z_l$-extension,
and let $\Delta$ be the order $p$ subgroup of
$\Gal(\tilde{\Q}/M)\simeq\Z/(l-1)$.
Here we note that $F\subseteq M$
since $l$ is unramified in $F/\Q$, which implies $\Delta\subseteq\mathcal{G}$.

Put $\mathcal{J}_{k_i}:=\Tor_{\Z_p}(X_{\tilde{k_i}}(p)_{\Delta})$.
Let $G:=\mathcal{G}/\Delta$, and let  
\[
l\mathcal{O}_{k_i}=\frak{l}_{i1}
\cdots\frak{l}_{ig_{k_i}(l)}\ \ \ (i=1,2)
\]
be the prime decomposition of $l$ in $k_i$ as above. 
Let $D_{ij}\subseteq \Gal(\tilde{k_i}/k_i)\simeq G$ be the decomposition subgroup for the prime $\frak{l}_{ij}$ $(i=1,2,\ 1\le j\le g_{k_i}(l))$.
Then we get the following exact sequence of $\F_p[[G]]$-modules by Proposition 3 and \eqref{ex}:
\begin{equation}\label{fe}
0\longrightarrow \mathcal{J}_{k_i}\longrightarrow\bigoplus_{j=1}^{g_{k_i}(l)}\F_p[[G/\pi_i(D_{ij})]]
\longrightarrow \F_p\longrightarrow 0.
\end{equation}
%
%

Now we require the following Lemma:
\begin{lem}
There exists an element $\sigma\in G$ 
such that 
the orders of $\sigma\pi_i(D_{ij})\in G/\pi_i(D_{ij})$ is finite and divisible by $p$ for any $1\le j\le g_{k_i}(l)$ and $i=1,2$.
\end{lem}

\

{\bf Proof.}\ \ \ 
Let $D_l\subseteq G=\Gal(\tilde{\Q}^{\Delta}/F)$ be the decomposition subgroup
for a prime of $F$ lying over $l$
(Note that $D_l$ does not depend on the choice of a prime of $F$ lying
over $l$). 
Then we find that $\pi_i(D_{ij})\subseteq D_l$ and 
$[D_l:\pi_i(D_{ij})]$ is finite since $[D_l:\pi_i(D_{ij})]$ equals the inertia degree of $\frak{l}_{ij}$
in $k_i/F$.
Therefore if we choose $\sigma\in G$ so that the order of $\sigma D_l\in
G/D_l$ is finite and divisible by $p$
(such $\sigma$ certainly exists),
then the order of $\sigma \pi_i(D_{ij})\in G/\pi_i(D_{ij})$ is also finite and 
divisible by $p$ for any $1\le j\le g_{k_i}(l)$ and $i=1,2$.
\hfill$\Box$

\

Let $\sigma\in G$ be an element given by Lemma 2.
It follows from \eqref{fe} and Lemma 1(2) that 
\[
\bigoplus_{j=1}^{g_{k_i}(l)}\F_p[[G/\pi_i(D_{ij})]]^{\langle\sigma\rangle}
\subseteq I_G\left(\bigoplus_{j=1}^{g_{k_i}(l)}
\F_p[[G/\pi_i(D_{ij})]]\right)\subseteq\mathcal{J}_{k_i}\ \ \ (i=1,2).
\]
Hence we have
\[
\mathcal{J}_{k_i}^{\langle\sigma\rangle}
=\bigoplus_{j=1}^{g_{k_i}(l)}\F_p[[G/\pi_i(D_{ij})]]^{\langle\sigma\rangle}\ \ \ (i=1,2).
\]
Therefore we obtain
\begin{equation}\label{jg}
(\mathcal{J}_{k_i}^{\langle\sigma\rangle})_G
\simeq\F_p^{\oplus g_{k_i}(l)}\ \ \ (i=1,2)
\end{equation}
by Lemma 1(3).

It follows from the assumption $X_{\tilde{k_1}}(p)\simeq X_{\tilde{k_2}}(p)$
as $\Z_p[[\mathcal{G}]]$-modules that $\mathcal{J}_{k_1}\simeq\mathcal{J}_{k_2}$
as $\F_p[[G]]$-modules.
Therefore we conclude that $g_{k_1}(l)=g_{k_2}(l)$ from \eqref{jg}.
\hfill$\Box$

\

Let $N/\Q$ be the Galois closure of $k_1k_2/\Q$.
Then $\Q(\mu_p)\cap N=\Q$ holds by the assumption of Theorem 1.
Hence For any given prime number $q\ne p$ unramified in $N$, there exists a prime number $l$
such that $[l,N/\Q]=[q,N/\Q]$, $l\equiv 1\pmod {p}$,
and $l$ is unramified in $k_1k_2$ by the \v{C}ebotarev density theorem, where $[r,N/\Q]$ stands for the Frobenius conjugacy class of $r$ in $\Gal(N/\Q)$ for any prime number $r$.
Then we have
\begin{equation}\label{g}
g_{k_1}(q)=g_{k_1}(l)=g_{k_2}(l)=g_{k_2}(q)
\end{equation}
for all but finitely many prime numbers $q$ by Proposition 4.

Finally, we employ the following theorem:
\begin{thm}[D. Stuart-R. Perlis \cite{S-P}]
For number fields $F_1$ and $F_2$ of finite degree,
$\zeta_{F_1}(s)=\zeta_{F_2}(s)$ holds if and only if
$g_{F_1}(q)=g_{F_2}(q)$ holds for all 
prime numbers $q$ with the possible exception of set primes of
Dirichlet density zero.
\end{thm}
Thus we conclude that $\zeta_{k_1}(s)=\zeta_{k_2}(s)$ holds by \eqref{g}
and the above theorem.\hfill$\Box$

\

We will close this section by giving following corollary to
Theorem 1:
\begin{cor}
Let $p$ be a prime number, and let 
$k_1$ and $k_2$ be number fields of finite degree.
We assume that $F:=k_1\cap\tilde{\Q}=k_2\cap\tilde{\Q}$ and that
$N\cap\Q(\mu_p)=\Q$, where $N$ is the Galois closure of $k_1k_2$ over $\Q$.
Furthermore we assume that $k_1/\Q$ is a Galois extension.
If $X_{\tilde{k_1}}(p)\simeq X_{\tilde{k_2}}(p)$ as $\Z_p[[\Gal(\tilde{\Q}/F)]]$-modules,
then $k_1=k_2$ holds.
\end{cor}
{\bf Proof.}\ \ \ 
It follows from $X_{\tilde{k_1}}(p)\simeq X_{\tilde{k_2}}(p)$ as $\Z_p[[\Gal(\tilde{\Q}/F)]]$-modules
that $\zeta_{k_1}(s)=\zeta_{k_2}(s)$ by Theorem 1.
Since $k_1/\Q$ is Galois extension,
the above equality implies $k_1=k_2$ by
\cite[p.93, Theorem (1.16)]{K}.\hfill$\Box$

\

One may expect that $X_{\tilde{k}}(p)$ determines the isomorphism class of $k$ by the above Corollary.
However we will see that this is not a case in general in the next section.

\section{Proof of Theorem 2}
In this last section, we will give a proof of Theorem 2.
We require the following:
\begin{lem}\label{lem-quot}
Let $k_1$ and $k_2$ be number fields of finite degree, 
and let $N/\Q$ be a finite Galois extension with $k_1k_2\subseteq N$.
Put $G:=\Gal(N/\Q),\ H_i:=\Gal(N/k_i)\ (i=1,2)$.
Assume that $\zeta_{k_1}(s)=\zeta_{k_2}(s)$ and that a prime number $p$ is prime
to $\#G$.
Then we have
\[
M_{H_1}\simeq M_{H_2} \mbox{ as $\Z_p[\mathcal{G}]$-modules}
\]
for any group $\mathcal{G}$ and $\Z_p[G\times\mathcal{G}]$-module $M$.
\end{lem}
{\bf Proof.}
It follows from the assumption $\zeta_{k_1}(s)=\zeta_{k_2}(s)$ that 
\[
\Q_p[G/H_1]\simeq\Q_p[G/H_2] \mbox{ as $\Q_p[G]$-modules}
\]
by \cite[P. 77, Theorem (1.3)]{K}, which in turn induces the isomorphism
\[
\varphi:\ \Z_p[G/H_1]\simeq\Z_p[G/H_2]\mbox{ as $\Z_p[G]$-modules}
\]
by \cite[p.626, (30.16)]{CR} because $p\nmid\# G$. 

Assume that $\varphi(1H_1)=\sum_{i=1}^{[G:H_2]}c_i(g_iH_2)$
($c_i\in\Z_p$) holds, and put
$\alpha:=\sum_{i=1}^{[G:H_2]}c_ig_i\in\Z_p[G]$,
where $G/H_2=\coprod_{i=1}^{[G:H_2]} g_iH_2$.

On the other hand, we find that
\begin{align*}
\psi_i:M_{H_i}&\simeq 
M\otimes_{\Z_p[G]}\Z_p[G/H_i]\\
x+\sum_{h\in H_i}(h-1)M&\mapsto x\otimes 1H_i
\end{align*}
as $\Z_p$-modules. 
Here, we define the right $\Z_p[G]$-module structure
on $M$ by $x\beta:=\beta^*x$, where
$\beta^*\in\Z_p[G]$ is defined to be $\sum_{\sigma\in G}c_\sigma\sigma^{-1}$ for
$\beta=\sum_{\sigma\in G}c_\sigma\sigma$.
Hence we get
\begin{equation}\label{iso}
M_{H_1}\overset{\psi_1}{\simeq} M\otimes_{\Z_p[G]}\Z_p[G/H_1]
\overset{\mathrm{id_M}\otimes\varphi}{\simeq} 
M\otimes_{\Z_p[G]}\Z_p[G/H_2]
\overset{\psi^{-1}_2}{\simeq}
M_{H_2}
\end{equation}
as $\Z_p$-modules, 
which maps $x\mod\sum_{h\in H_1}(h-1)M\in M_{H_1}$
to $\alpha^*(x)\mod\sum_{h\in H_2}(h-1)M\in M_{H_2}$.

Because each element of $G$ and that of $\mathcal{G}$ commute,
the $\Z_p$-module isomorphism given by \eqref{iso}
is in fact a $\Z_p[\mathcal{G}]$-isomorphism.
\hfill$\Box$

\

{\bf Proof of Theorem 2}\ \ \ 
We first note that the restriction induces the isomorphism
\[
\Gal(\tilde N/\Q)\simeq\Gal(\tilde{\Q}/\Q)\times\Gal(N/\Q),
\]
since $N\cap\tilde{\Q}=\Q$.
Then $\Gal(\tilde{\Q}/\Q)\times\Gal(N/\Q)$
operates on $X_{\tilde{N}}(p)$.
Let $H_i:=\Gal(N/k_i)\ (i=1,2)$.
It follows from Lemma 3 that there exists $\Z_p[[\Gal(\tilde{\Q}/\Q)]]$-module isomorphism
\begin{equation}\label{H_i}
X_{\tilde{N}}(p)_{H_1}\simeq X_{\tilde{N}}(p)_{H_2}.
\end{equation}

We find that $\tilde{N}^{H_i}=\tilde{k_i}$
and $X_{\tilde{N}}(p)_{H_i}=X_{\tilde{N}}(p)_{\Gal(\tilde{N}/\tilde{k_i})}$.
Because $p\nmid \#\Gal(\tilde{N}/\tilde{k_i})$, we see that
\begin{equation}\label{co}
X_{\tilde{N}}(p)_{\Gal(\tilde{N}/\tilde{k_i})}\simeq X_{\tilde{k_i}}(p)\ \ \ 
\end{equation}
as $\Z_p[[\Gal(\tilde{k_i}/k_i)]]\simeq\Z_p[[\Gal(\tilde{\Q}/\Q)]]$-modules $(i=1,2)$.
Thus \eqref{H_i} and \eqref{co} imply the assertion of the theorem.
\hfill$\Box$

\

{\bf Example} (see \cite[p.359]{P})\ \ \ 
Define irreducible polynomials
\[
f_1(X)=X^7-7X+3,\ f_2(X)=X^7+14X^4-42X^2-21X+9.
\]
over $\Q$, and let $\alpha_i\in\Qb$ be a root of $f_i(X)\ (i=1,2)$.
Put $k_i:=\Q(\alpha_i)\ (i=1,2)$.
Then $k_1$ and $k_2$ are number fields of degree 7 and
they are not isomorphic.
Furthermore $\zeta_{k_1}(s)=\zeta_{k_2}(s)$ holds. 
Their common Galois closure $N$ over $\Q$ is a 
$GL_3(\F_2)$-extension over $\Q$.
Because $\Gal(N/\Q)\simeq GL_3(\F_2)$ is a finite simple group of order $168=2^3\cdot 3\cdot 7$,
we see $N\cap\tilde{\Q}=\Q$.
Hence it follows from Theorem 2 that 
\[
X_{\tilde{k_1}}(p)\simeq X_{\tilde{k_2}}(p)
\]
as $\Z_p[[\Gal(\tilde{\Q}/\Q)]]$-modules
for any prime number $p$ different from $2,3$ and $7$.

The above example says that $X_{\tilde{k}}(p)$ for a``good" prime $p$
determines $\zeta_k(s)$ (Theorem 1), but can not determine the isomorphism class of $k$.

\

\

\noindent
Manabu Ozaki,\par\noindent
Department of Mathematics,\par\noindent
School of Fundamental Science and Engineering,\par\noindent
Waseda University,\par\noindent
Ohkubo 3-4-1, Shinjuku-ku, Tokyo, 169-8555, Japan\par\noindent
e-mail:\ \verb+ozaki@waseda.jp+
\end{document}